\documentclass{amsart}
\usepackage{amsfonts}
\pdfoutput=1
\usepackage{graphicx}
\usepackage{tabularx}
\usepackage{array}
\usepackage[usenames,dvipsnames]{color}
\usepackage{comment}
\usepackage{xcolor,soul}
\usepackage {tikz}
\usetikzlibrary {positioning}

\begin{document}
\newcommand{\cbl}{\color{blue}}
\newcommand{\cro}{\color{red}}
\newcommand{\cno}{\color{black}}

\newcommand{\C}{\mathbb{C}}
\newcommand{\F}{\mathbb{F}}
\newcommand{\R}{\mathbb{R}}
\newcommand{\Q}{\mathbb{Q}}
\newcommand{\N}{\mathbb{N}}
\newcommand{\Z}{\mathbb{Z}}
\newcommand{\PP}{\mathbb{P}}
\newcommand{\HH}{\mathbb{H}}
\newcommand{\MU}{\mathbb{\mu}}
\newcommand{\TT}{\mathbb{T}}
\newcommand{\GG}{\mathbb{G}}
\newcommand{\SSS}{\mathbb{S}}
\newcommand{\m}{\mathrm{m}}

\newtheorem{thmA}{Theorem}
\newtheorem{defnA}[thmA]{Definition}
\newtheorem{propA}[thmA]{Proposition}
\newtheorem{coroA}[thmA]{Corollary}
\newtheorem{lemA}[thmA]{Lemma}
\newtheorem{conjA}[thmA]{Conjecture}
\newtheorem{rem}[thmA]{Remark}
\newtheorem{ex}[thmA]{Example}
\newtheorem{exs}[thmA]{Examples}
\newtheorem{pbmA}[thmA]{Problem}
\newcommand{\dd}{\;\mathrm{d}}
\newcommand{\e}{\mathrm{e}}
\newcommand{\pf}{\noindent {\bf PROOF.} \quad}
\newcommand{\res}{\mathrm{Res}}

\newcommand{\Li}{\mathrm{Li}} 
\newcommand{\I}{\mathrm{I}} 
\newcommand{\re}{\mathop{\mathrm{Re}}} 
\newcommand{\im}{\mathop{\mathrm{Im}}} 
\newcommand{\ii}{\mathrm{i}} 
\newcommand{\Lf}{\mathrm{L}} 
\newtheorem{obs}[thmA]{Observation} 

\newcommand{\frank}[1]{{\color{blue}{Frank: #1}} }
\newcommand{\matilde}[1]{{\color{red}{Matilde: #1}} }

\title{Edge clique covers in graphs with independence number two: a special case}

\author{Frank Ramamonjisoa}
\address{D\'{e}partment de Math\'{e}matique et de Statistique, Universit\'{e} de Montr\'{e}al. CP 6128, succ. Centre-Ville. Montreal, QC H3C 3J7, Canada}
\email{frank.ramamonjisoa@umontreal.ca}

\thanks{The authors of this work have been supported by the Natural Sciences and Engineering Research Council of Canada (BESC D)}

\keywords{edge clique cover number, intersection number, claw-free graphs}

\begin{abstract} The edge clique cover number $ecc(G)$ of a graph $G$
is the size of the smallest set of complete subgraphs whose union covers
all edges of $G$. It has been conjectured that all the simple graphs
with independence number two satisfy $ecc(G)\leq n$. First, we present
a class of graphs containing edges
difficult to cover but that satisfy the conjecture. Second, we
describe a large class of graphs $G$ such that
$ecc(G)\leq \frac{3}{2}n$. This class is easy to
characterize.  \end{abstract}

\maketitle

\section{Introduction} 

In this paper, $G$ will be a simple graph and $n$ will be the number
of vertices of $G$. A \emph{clique} of $G$ is a subset of the
vertices of $G$ that induces a complete graph.
The \emph{edge clique cover number} of $G$, denoted $ecc(G)$, is the
minimum number of cliques in $G$ such that both endpoints of each
edge of $G$ lie in at least one of them.  At the beginning was the
following problem:
\begin{pbmA} (Chen, Jacobson, Kézdy, Lehel, Scheinerman, Wan \cite{2}) If $G$ is a claw free graph, is $ecc(G)\leq n$?
\end{pbmA}
A {\em claw} is a graph with four vertices and three edges linking
three of the vertices to the fourth. Javadi and Hajebi gave a positive
answer to the question for graphs with independence number greater than two
\cite{6}. Their proof is based on Chudnovsky and Seymour’s
structure theorem for claw-free graphs (see \cite{3}), concerning
only graphs with independence number $\alpha(G)\geq 3$. This leaves  the following conjecture:     
\begin{conjA}\label{two} If $G$ is a simple graph with $\alpha(G)=2$, then $ecc(G)\leq n$.
\end{conjA}

Special classes of graphs satisfying this conjecture have been found.
For example graphs with diameter three, or $2K_2$-free and co-claw
free graphs \cite{1}. Our first result adds one more class of graphs
that satisfy the conjecture. It contains  graphs with dominating edges, which
are particularly difficult to cover in graphs with independence
number two and diameter two.     
\begin{defnA}
An edge $uv\in E(G)$ is a \emph{dominating edge} of $G$ if every
other vertex of $G$ is adjacent to at least one of $u$ and $v$.
\end{defnA}
We can now state our first theorem:
 \begin{thmA} \label{premier} Let $G$ be a simple graph with $\alpha(G)=2$. If $G$ has a
dominating edge $uv$ such that $G\setminus \{u,v\}$ has diameter 3, then
$ecc(G)\leq n$.  
\end{thmA}
Thus to prove the conjecture we only have to do so for the class of
graphs $G$ with independence number two, diameter two, and such that
for any 
dominating edge $uv$ the graph
$G\setminus \{u,v\}$ has diameter two.

As this is still out of reach, we can try to find upper bounds for the
edge clique cover number. The best known bound on $ecc(G)$
for general graphs with independence number two is the minimum of
$n+\delta(G)$ and $2n-\Omega(\sqrt{n\log n})$, where $\delta(G)$ is
the minimum degree of $G$ \cite{1}. Better bounds has been found for some classes of graphs. For example, it has been shown in
\cite{1} that if the chromatic number $\chi(\bar{G})$ of the complement of $G$ is $4$ then $ecc(G)\leq \frac{3}{2}n+4$.
The combination of a large class of graphs $G$ with a good bound on
the  $ecc(G)$ is mostly elusive at this time. \\

Our second result is
about one such combination with the additional pleasant property that
it is easy to describe: avoid incompatible
double top and supercycle free graphs.
\begin{defnA}
Let $G$ be a simple graph and $x_0,x_1,x_2,x_3,x_4$ be vertices of
$G$. The  path
$x_0x_1x_2x_3x_4$ is a {\bf double top} if its edges are dominating
and $x_2$ is adjacent neither to $x_0$ nor
to $x_4$ in $G$.
\end{defnA} 

Examples of double tops. The dominating edges are in red.
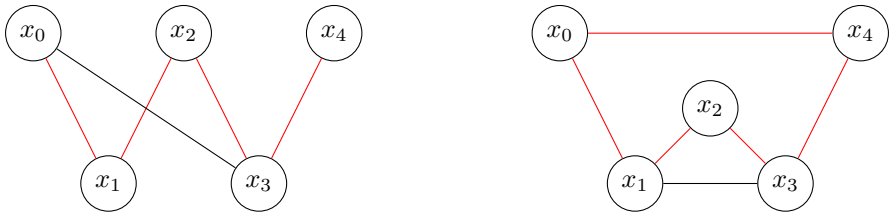
\begin {center}
\begin {tikzpicture}
\node[draw, circle] (C1) at (0,0) {$x_0$}; \node[draw, circle] (C2) at (1,-2) {$x_1$}; \node[draw, circle] (C3) at (2,0) {$x_2$}; \node[draw, circle] (C4) at (3,-2) {$x_3$}; \node[draw, circle] (C5) at (4,0) {$x_4$}; 
\draw[red] (C1) -- (C2) ; \draw[red]  (C2) to (C3) ;\draw[red]  (C3) to (C4) ;\draw[red]  (C4) to (C5) ;\draw[]  (C1) to (C4) ;

\node[draw, circle] (D1) at (7,0) {$x_0$}; \node[draw, circle] (D2) at (8,-2) {$x_1$}; \node[draw, circle] (D3) at (9,-1) {$x_2$}; \node[draw, circle] (D4) at (10,-2) {$x_3$}; \node[draw, circle] (D5) at (11,0) {$x_4$}; 
\draw[red] (D1) -- (D2) ; \draw[red]  (D2) to (D3) ;\draw[red]  (D3) to (D4) ;\draw[red]  (D4) to (D5) ;\draw[]  (D2) to (D4) ; \draw[red]  (D1) to (D5) ; 

\end{tikzpicture}
\end{center}

If $G$ has few  dominating edges, double tops should not be too common
and avoiding graphs with double tops should not be  too big a
constraint. In fact, we can  relax the constraint somewhat.
\begin{defnA}
A double top $x_0x_1x_2x_3x_4$ is \emph{incompatible} if $x_1x_3$ is not an edge.
\end{defnA}  

Example of incompatible double top. The dominating edges are in red.
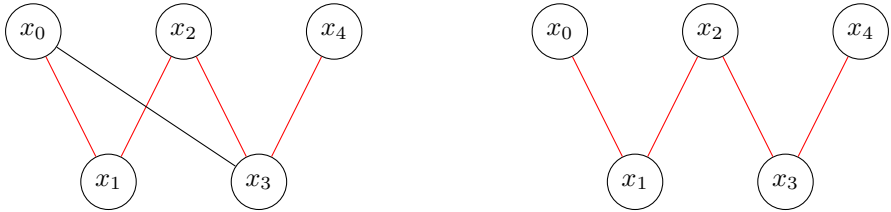
\begin {center}
\begin {tikzpicture}
\node[draw, circle] (C1) at (0,0) {$x_0$}; \node[draw, circle] (C2) at (1,-2) {$x_1$}; \node[draw, circle] (C3) at (2,0) {$x_2$}; \node[draw, circle] (C4) at (3,-2) {$x_3$}; \node[draw, circle] (C5) at (4,0) {$x_4$}; 
\draw[red] (C1) -- (C2) ; \draw[red]  (C2) to (C3) ;\draw[red]  (C3) to (C4) ;\draw[red]  (C4) to (C5) ;\draw[]  (C1) to (C4) ;

\node[draw, circle] (C1) at (7,0) {$x_0$}; \node[draw, circle] (C2) at (8,-2) {$x_1$}; \node[draw, circle] (C3) at (9,0) {$x_2$}; \node[draw, circle] (C4) at (10,-2) {$x_3$}; \node[draw, circle] (C5) at (11,0) {$x_4$}; 
\draw[red] (C1) -- (C2) ; \draw[red]  (C2) to (C3) ;\draw[red]  (C3) to (C4) ;\draw[red]  (C4) to (C5) ;
\end{tikzpicture}
\end{center}
We will then allow double tops in general but avoid incompatible ones.
We need one more definition.

\begin{defnA}
We call {\bf supercycle} a cycle $u_0\cdots u_k,~k\in \{3,\cdots,
n-1\}$ of dominating edges such that for all $i\in \{0,\cdots, k\}$,
$u_{i-1}u_{i+1}$ is an edge, addition modulo $k+1$.
\end{defnA}

Example of supercycle.
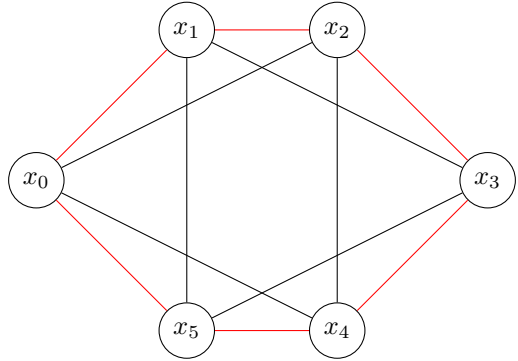
\begin {center}
\begin {tikzpicture}

\node[draw, circle] (C1) at (0,0) {$x_0$}; \node[draw, circle] (C2) at (2,2) {$x_1$}; \node[draw, circle] (C3) at (4,2) {$x_2$}; \node[draw, circle] (C4) at (6,0) {$x_3$}; \node[draw, circle] (C5) at (4,-2) {$x_4$}; \node[draw, circle] (C6) at (2,-2) {$x_5$};
\draw[red] (C1) -- (C2) ; \draw[red]  (C2) to (C3) ;\draw[red]  (C3) to (C4) ;\draw[red]  (C4) to (C5) ;\draw[red]  (C5) to (C6) ; \draw[red]  (C6) to (C1) ;
\draw[]  (C1) to (C3) ;\draw[]  (C2) to (C4) ;\draw[]  (C3) to (C5) ;\draw[]  (C4) to (C6) ;\draw[]  (C5) to (C1) ; \draw[]  (C6) to (C2) ;
\end{tikzpicture}
\end{center}

So a supercycle is a cycle $C$ of at least four  dominating
edges in which any two vertices at distance two on $C$ form an edge. We are
now ready to state the following result:
\begin{thmA}\label{deuxieme}
If $G$ is a simple graph with $\alpha(G)=2$ containing neither supercycles nor incompatible double tops, then $ecc(G)\leq \frac{3}{2}n $.
\end{thmA}

This document is organized as follows. In Section 2 we consider some
preliminary and/or known results on the edge clique covering number of
graphs with independence number two. In Section 3 we prove
Theorem~\ref{premier}. In Section 4 we prove
Theorem~\ref{deuxieme} and we conclude in Section 5. 

\section{Typology of graphs with $\alpha(G)=2$}

Conjecture \ref{two} concerns graphs with independence number two. Let
us first look at some known results on the edge clique covering number of
such graphs.\\

If a graph $G$ is not connected, then it cannot have more
than two connected components because $\alpha(G)=2$, and every connected component is a
clique. Thus  no generality is lost in assuming that $G$ is
connected.\\

A graph $G$ with $\alpha(G)=2$ has diameter $2$ or $3$. It can
clearly not be  of diameter $1$ and if it had diameter  $4$ or more,
it would have $\alpha(G)\geq 3$. The diameter $3$ case has been solved
in \cite{1}.
\begin{propA}[\cite{1}]
Let $G=(V,E)$ be a simple graph, if $\alpha(G)=2$ and \\$diam(G)=3$,
then $ecc(G)\leq \lceil \frac{n+1}{2}\rceil$.
\end{propA}
This leaves only graphs of diameter $2$ (and of course
$\alpha(G)=2$) to consider. Since the conjecture is trivially true for $n\leq 3$, we will suppose that $n\geq 4$.\\ 
As usual, let us denote by $\overline{N[u]}$ the set of non-neighbours of $u$. In a graph with independence
number two, each $\overline{N[u]}$ is a clique, possibly trivial. 
If there is no dominating edge in $G$, the set
$\{\overline{N[u]}\}_{u\in V(G)}$ is an edge clique cover of $G$ that
contains at most $n$ cliques. Indeed, for each edge $xy$ of $G$, there
is a vertex $u$ that is not adjacent to either of its ends, and so $xy$
is covered by the clique $\overline{N[u]}$. We thus have the following
lemma, proved in \cite{1}:
\begin{lemA}[\cite{1}]
Let $G$ be a simple graph such that $\alpha(G)=2$. If $G$ has no dominating edge, then $ecc(G)\leq n$.
\end{lemA}

We will suppose from now on that $G$ possesses a dominating edge
$uv$. Thus, each vertex $w$ of $G$ (different from $u$ and $v$) is adjacent to $u$ or adjacent to $v$ (or adjacent to both). Let $U$ be the set of vertices adjacent to $u$ but not adjacent to $v$, $V$ be the set of vertices adjacent to $v$ but not adjacent to $u$, and $W$ be the set of vertices adjacent to both $u$ \emph{and} $v$. Without loss of generality,
we will suppose that $|U|\geq |V|$. Note that both $U$ and $V$ are
cliques.\\

Let's consider the subgraph $G'=G \setminus
\{u,v\}$. We have
$\alpha(G')\leq 2$ because $\alpha(G)= 2$. If $\alpha(G')=1$, then $V(G')$ is
a clique and the edges of $G$ are covered by the four cliques $V(G'), \{u,v\}, U\cup W \cup
\{u\}, V\cup W \cup \{v\}$ and the conjecture is true.

\begin{lemA}
Let $G$ be a simple graph with $\alpha(G)=2$. If $G$ has a dominating
edge $uv$ such that $G\setminus\{u,v\}$ is of diameter 1, then $ecc(G)\leq n$.
\end{lemA}

If $\alpha(G')=2$, suppose first that $G'$ is not connected. As we have seen before, there can only be two
connected components, $C_1$ and $C_2$, and those connected components
are cliques. As $U$ and $V$ are cliques, each of them is included in
one connected component. Suppose for example that $U
\subset C_1$ and $V \subset C_2$. Then, $G$ is covered by the cliques:
$V(C_1)\cup \{u\}, V(C_2)\cup \{v\}, \{W\cap V(C_1)\}\cup \{v\}, \{ W\cap
V(C_2)\}\cup \{u\}, \{u,v\}$. \\ As $G$ is covered by $5$ cliques, the conjecture holds when $n\geq 5$. If $U$ and $V$ are included in the same connected component, this is still true. Finally, we can check directly that if $n=4$, this is also true. So we will suppose from now on that $G'$ is
connected.\\

We are now facing two cases, the diameter of $G'$ is 3 or the
diameter of $G'$ is 2, and the purpose of the next section is to
study the first one.\\

\section{Proof of Theorem \ref{premier}}
When the diameter of $G'$ is three, then according to Proposition $9$,
the edges of $G'$ are covered by at most $\lceil \frac{n(G')+1}{2}
\rceil=\lceil \frac{n-1}{2} \rceil \leq \frac{n}{2}$ cliques $\{D^i\}_{1\leq i\leq k},
~k\in\mathbb{N}$, $n(G')$ being the number of vertices of the graph $G'$.\\

Let's consider the graph $G[W]$ induced by the vertices of $W$. We
need to cover the vertices of $W$ by cliques of $G[W]$ in
the rest of the proof, in order to cover the edges between
$u,v$ and $W$. We have a bound for the number of those
cliques: 
\begin{lemA} The vertices of $W$ can be covered by cliques $\{C^i\}_{1\leq i\leq
\ell},~\ell\in\mathbb{N}$ of $G[W]$, with $l\leq \lceil
\frac{|W|+1}{2} \rceil$.  
\end{lemA}

\begin{proof}
Suppose first that $G[W]$ is not connected. As we've seen before, in that
case, $G[W]$ consists of two cliques and the lemma is true. So we will suppose that $G[W]$ is connected.\\ 
We will also suppose that $\alpha(G[W])=2$ because if $\alpha(G[W])=1$ or if $W$ is empty, the lemma is obviously satisfied.\\
If $G[W]$ has only one or two vertices, the lemma is true.\\
If $G[W]$ has only 3 vertices, there are two adjacent vertices and the lemma is true.\\
If $G[W]$ has only 4 vertices, then we can check that we can cover them with two sets and the lemma is true.\\
\indent

Suppose now that $G[W]$ has more than 4 vertices. Let's take any
vertex $w$ of $W$. It lies in a clique of
$G[W]$ that has at least two vertices since $G[W]$
is connected. We call this clique $C^1$ and we consider the graph $G^1=G[W]\setminus\{C^1\}$.\\
If $G^1$ is not connected, $G^1$ satisfies the
lemma and thus, $G[W]$ satisfies the lemma. So let's suppose that
$G^1$ is connected. If $G^1$ has more than 4 vertices, then we repeat the above process by
taking two of the remaining vertices linked by an edge and by
considering a clique $C^2$ that contains them. As we remove at each
step at least two vertices, if the remaining graph possesses at a
certain step 3 or 4 vertices, the lemma is true. Otherwise we go from
a graph with more than 4 vertices to a graph with at most two
vertices, and in that case, the lemma is also satisfied, because we
would have removed at least 3 vertices in one step before arriving to
at most two vertices. \\
\indent
In conclusion, the lemma holds in all cases.

\end{proof}
Let $\{C^i\},\ 1\leq i \leq \ell,\ \ell \in \mathbb{N} $ be the
set of cliques that cover the vertices of $G[W]$ found in the preceeding Lemma. For
each $i$, let $B^i=C^i\cup \{u,v\}$. The set $\{B^i: 1\leq i \leq
\ell  \in \mathbb{N}\}$
 covers $G[W]$ and all the edges between
$u,v$ and $W$. If $W$ is not empty, this set of cliques is not
empty, which means that it also covers the edge $uv$, since $G[W]$ is not connected and $\alpha(G)=2$. If $W$ is empty,
we consider the set composed of only the clique
$\{u,v\}$.\\ If $U$ is 
not empty, $U'=U\cup \{u\}$ is a clique and it covers all the edges
from $u$ to $U$; if $U$ is empty, we set $U'=\emptyset$. We do the
same for $V$.\\
The cliques $\{D^i\}_{1\leq i\leq k}\cup \{B^i\}_{1\leq i\leq \ell} \cup \{U'\} \cup \{V'\}$ cover all the edges of $G$ and there are at most 
\[m=\lceil \frac{n-1}{2} \rceil+\lceil \frac{|W|+1}{2} \rceil+1_{U'\neq \emptyset}+1_{V'\neq \emptyset}\hspace*{4cm} (1)\]
of them. The notation $1_{P}$ refers to the indicator function on the assertion $P$: this function is equal to $1$ if and only if $P$ is true.\\

 Let's study this formula with different values of $1_{U'\neq
 \emptyset}$ and $1_{V'\neq \emptyset}$.

\begin{itemize}
\item If $U$ and $V$ are empty, $m=\lceil \frac{n-1}{2} \rceil+\lceil \frac{n-2+1}{2} \rceil\leq \frac{n}{2}+\frac{n}{2}\leq n$. We have thus $ecc(G)\leq n$.
\item If only one of the sets $U$ and $V$ is not empty (suppose that it's $U$ without loss of generality), as $|W|\leq n-3$, we have:\\
 $m\leq \lceil \frac{n-1}{2} \rceil+\lceil \frac{n-3+1}{2} \rceil+1\leq \lceil \frac{n-1}{2}\rceil+\lceil \frac{n-2}{2}\rceil+1\leq \frac{n}{2}+\frac{n-1}{2}+1\leq n+\frac{1}{2}.$\\
But as $ecc(G)$ is an integer, we have also in this case $ecc(G)\leq n$.
\item If $U$ and $V$ are both non empty, we have then $|W|\leq n-4$. Suppose that $U\cup V$ contains at least 3 vertices. We have then $|W|\leq n-5$ and:\\
 $m\leq \lceil \frac{n-1}{2}\rceil+\lceil \frac{n-5+1}{2}\rceil+2\leq \frac{n}{2}+\frac{n-3}{2}+2\leq n+\frac{1}{2}$. As $ecc(G)$ is an integer, we deduce $ecc(G)\leq n$.
\end{itemize}

 So now there is only one subcase to investigate in order to finish the proof: $|U\cup V|=2$ and $U$ and $V$ both non empty, so $|U|=|V|=1$. We will first improve Lemma 12:
\begin{lemA} If $|W|\geq 4$ and if $G[W]$ is not a 5-cycle, the vertices of $W$ can be covered by at most $\lceil \frac{|W|}{2} \rceil$ sets $\{C^i\}_{1\leq i\leq \ell},~\ell\in\mathbb{N}$, a set $C^i$ being a clique of $G[W]$.
\end{lemA}
\begin{proof}
We have seen before that if $G[W]$ is not connected, it consists
of two cliques and the current lemma is true. We
will then suppose, without loss of generality, that $G[W]$ is connected.\\
As in Lemma 12, we see that if $G[W]$ has four vertices, we can cover
them by at most two sets, and the current lemma is true. \\
We can check directly that if $|W|=5$ and $G[W]$ is not a 5-cycle,
$G[W]$ contains a clique of at least 3 vertices.\\
 Suppose now that $|W|> 5$. There exists a connected subgraph $G''$ of
$G[W]$ with exactly five vertices because $G[W]$ is connected. If $G''$ is a 5-cycle, then since $G[W]$ is connected, there exists a
connected subgraph of $G[W]$ containing $G''$ with exactly six
vertices. We can check directly that this subgraph possesses a clique
 of at least 3 vertices. If $G''$ is not a
5-cycle, we have seen above that it too contains a clique  of
at least 3 vertices. In either case, $G''$ contains a clique
$\mathcal{C}$ of at least three vertices.\\
  From Lemma 12 we have that the vertices of $G[W]\setminus\mathcal{C}$ are covered by at most $\lceil \frac{|V(G[W]\setminus C)|+1}{2} \rceil\leq \lceil \frac{|W|-3+1}{2} \rceil$ sets $\{C^i\}_{1\leq i\leq \ell},~\ell\in\mathbb{N}$. So the vertices of $G[W]$ are covered by at most $\lceil \frac{|W|-3+1}{2} \rceil+1=\lceil \frac{|W|-3+3}{2} \rceil$ sets $\{C^i\}_{1\leq i\leq \ell}\cup C$.
\end{proof}
So if $|W|\geq 4$ and if $G[W]$ is not a 5-cycle, we can adapt formula
$(1)$ to our present case and we have at most \\
$m=\lceil \frac{n-1}{2}
\rceil+\lceil \frac{|W|}{2} \rceil+1+1= \lceil \frac{n-1}{2}
\rceil+\lceil \frac{n-4}{2} \rceil+2\leq
\frac{n}{2}+\frac{n-3}{2}+2\leq n+\frac{1}{2}$ \\cliques covering the
edges of $G$. As the number of cliques is an integer, the theorem
is true.\\
Now if $|W|= 3$, the vertices of $W$ are covered by at most two sets. As in addition $n=7$, formula $(1)$ gives $m=7$ and the theorem is true. If $|W|=2$ and it is covered by one clique, as $n=6$, formula $(1)$ gives $m=6$ and the theorem is true. If $|W|=1$, as $n=5$, formula $(1)$ gives $m=5$ and the theorem is true. And finally if $W$ is empty, the theorem is trivially satisfied.\\
\indent
If $W$ is composed of two non-adjacent vertices, we can check directly that all the different cases satisfy the theorem. And finally, if $G[W]$ is a 5-cycle, and by considering that the vertices in $G[W]$ not adjacent to the vertex in $U$ (or $V$) forms a clique, we also check directly in all the cases that the theorem is true. This concludes the proof of Theorem \ref{premier}.

\section{Proof of Theorem \ref{deuxieme}}

We saw in Section $2$ that if there is no dominating edge, the set of
cliques $D=\{\overline{N[x]}\}_{x\in V(G)}$ covers all the edges of
$G$. In a graph containing dominating edges, $D$ covers all the
non-dominating edges of $G$. But it covers no dominating edge. The
idea is to extend these cliques in order to cover as many dominating
edges as possible, and to count the additional cliques necessary to
cover the remaining dominating edges.\\

A vertex $y$ adjacent to a vertex $x$ can be used to
extend the clique $\overline{N[x]}$ if all the vertices of
$\overline{N[x]}$ are linked to $y$. This
gives the following.
\begin{lemA}
Let $x$ be a vertex in $G$. The clique $\overline{N[x]}$ can be extended by the vertex $y$ if and only if $xy$ is a dominating edge.
\end{lemA}
In this case, all the dominating edges between $y$ and
$\overline{N[x]}$ will be covered by this extended clique.\\
\indent But extending a clique like $\overline{N[x]}$ is not always
relevant, because it may not cover a dominating edge. Let us then first study the methods we have to extend a clique
$\overline{N[x]}$ so it covers a dominating edge.\\
\indent If $xy$ is a dominating edge and neither $x$ nor $y$ belongs
to another dominating edge, $xy$ cannot be covered by an extension of
a clique $\overline{N[z]}$, for any vertex $z$ of $G$. Indeed, this would mean that either $zx$ or $zy$ is a dominating
edge. It means that we can cover $xy$ by an extension of cliques of
$D$ only if there exists a vertex $z$ in $G$ such that $xz$ or $yz$ is
dominating.\\
\indent Let's suppose without loss of generality that $xz$ is
dominating. Recall that $xy$ is dominating.  
\begin{itemize} 
\item if
$y$ and $z$ are not adjacent, $z$ belongs to $\overline{N[y]}$ and
this clique can be extended to $x$ as $xy$ is dominating. In that
case, the extended clique covers $xz$. In the same way, $xy$
is covered by a clique extended from $\overline{N[z]}$.  
\item if $yz$
is an edge but not dominating, then as $z$ doesn't belong to
$\overline{N[y]}$, $xz$ cannot be covered by a clique extended from
$\overline{N[y]}$, because $\overline{N[y]}$ cannot be extended to
both $x$ and $z$.  
\item if $yz$ is a dominating edge, it can be
covered by extending $\overline{N[x]}$ to $y$ and $z$.\\
 \end{itemize}
So for a dominating edge $xy$, there are two cases where it can be covered by an extension of a clique $\overline{N[z]}$ of $D$ for a certain vertex $z$ in $G$ so that $xz$ is dominating. And it cannot be covered by one of these extensions in the other cases:
\begin{lemA}
A dominating edge $xy$ is covered by an extension of the clique \\$\overline{N[z]},~z \in V(G)$ if and only if:
\begin{itemize}
\item $zx$ is dominating and $z$ is not adjacent to $y$,
\item $zy$ is dominating and $z$ is not adjacent to $x$,
\item $zx$ and $zy$ are both dominating edges.
\end{itemize}
\end{lemA}
But maybe we can extend $\overline{N[z]}$ in several ways, to several vertices adjacent to $z$, in order to cover several dominating edges. For example, let's see what happens if we can extend a clique $\overline{N[z]},~z\in V(G)$ to a vertex $y_1$ and cover at the same time a dominating edge, say $y_1w_1$with $w_1 \in \overline{N[z]}$, and if there is another possibility of extension of $\overline{N[z]}$ to a vertex $y_2$, that covers also a dominating edge, say $y_2w_2$ with $w_2 \in \overline{N[z]}$. If $y_1y_2$ is not an edge, we cannot extend $\overline{N[z]}$ to both $y_1$ and $y_2$ and we have to choose only one of theses possible extensions.\\
\indent
By choosing graphs with incompatible double tops, we avoid the situation where $y_1y_2$ is not an edge. Indeed, in the above situation, $w_1y_1zy_2w_2$ forms a double top. As it is not incompatible, $y_1y_2$ is an edge. So we can extend $\overline{N[z]}$ to both $y_1$ and $y_2$.\\

Following this idea, for each vertex $x$ in $G$, we extend the
clique $\overline{N[x]}$ to a larger clique, denoted by $C^x$, by adding all the vertices $y$ such that $xy$ is dominating and there is a dominating edge incident to $y$ and to a vertex in $\overline{N[x]}$. If there is no such vertex, then $C^x=\overline{N[x]}$. \\
 The set of cliques $C=\{C^x\}_{x\in V(G)}$ contains $n$ elements that
cover a certain number of dominating edges in $G$. But there may still
exist dominating edges that are not covered by $C$. We now have to
find a way to count them to determine a bound for the edge clique
cover number of $G$. Let $D'$ be a set of cliques, initially empty.
It will receive the cliques that we need in order to cover the
dominating edges that are not covered by $C$.\\

Suppose that $zw_1$ is a not covered (dominating) edge. Then all the
other dominating edges $zw$, if any, are such that $w$ is adjacent to
$w_1$. This is in particular the case for all the dominating edges
$zw$ \emph{yet not covered}. Indeed, if it were not the case, then we
could cover $zw_1$ with $C^w$. Thus, the vertices of all the uncovered dominating edges
incident to $z$ form a clique, and we will denote this clique by
$C_z=\{z,w_1,\cdots,w_k\},~k\in\mathbb{N}$. Note that $C_z$ is not the
same as the extended clique $C^z$. \\

If there is a $w_i,~1\leq i\leq k$ in $C_z$ such that there is no
not covered edge incident to $w_i$ with the other end outside $C_z$,
we put $C_z$ in $D'$ and we label it with $z$ and $w_i$.
\begin{defnA} We call a vertex $z$ in $G$ a {\em favourable} vertex if there exists a vertex $w\in C_z$ such that each dominating edge of $G$ containing $w$ has its second vertex in $C_z$. We call also $w$ a {\em twin vertex} of $z$.
\end{defnA}

The clique $C_z$ covers all the edges not yet covered and incident to $z$ and those incident to $w_i$.
 Thus all the remaining not covered edges are inside a subgraph of $G$ with at most $n-2$ vertices. \\
 \indent
 We continue the same process for the remaining edges not covered. If there is a favourable vertex, we put the corresponding clique in $D'$ and we label it with this vertex and one of its twin vertices. The process ends when there is no favourable vertex.\\
 \indent
 Observe now that a vertex can appear in only one clique of $D'$. This property will give us a bound for $|D'|$, and consequently for $ecc(G)$.\\
 \indent
 We remark that at the beginning of the process of putting cliques in $D'$, a given vertex $z'$ may not be favourable, but that after a certain number of steps, it may become favourable because all the not covered edges incident to its twin vertex with the other end outside $C_{z'}$ are covered by cliques in $D'$.\\

Among the (possibly) remaining not covered edges, there is a cycle. Indeed, if we consider a connected component of the subgraph formed by these not covered edges, and if there is no cycle, this is a tree, and in that case, there is a vertex of degree one, which is a contradiction because this vertex should be in a clique of $D'$. We can even say that ``there are only cycles" because there is no vertex of degree one in this subgraph. Moreover, if there is a cycle $u_0\cdots u_k,~k\in \{2,\cdots, n-1\}$, then for all $i\in \{0,\cdots, k\}$, $u_{i-1}iu_{i+1}$ is an edge in $G$, addition modulo $k+1$. This is because all the remaining edges incident to a given vertex of $G$ form a clique. So we are in the presence of cycles of length $3$ or supercycles (we recall that a supercycle has length at least 4). \\

We take into account the previous stages: the cliques in $D$ are extended and we have put cliques in $D'$ until there is no remaining favourable vertex. If there remains uncovered edges, we have just seen that there is a cycle, and this cycle must be a supercycle or a cycle of length 3. But then it is a cycle of length three because by hypothesis there is no supercycle in $G$. We denote this cycle by $u_0u_1u_2$. As $u_1$ cannot be a favourable vertex, there is a vertex $u_3$ so that $u_2u_3$ is a dominating edge and $u_3$ doesn't belong to $C_{u_1}$. So we have a path in $G$ with 4 vertices. Again as $u_2$ is not a favourable vertex, there is a vertex $u_4$ so that $u_3u_4$ is a dominating edge and $u_4$ doesn't belong to $C_{u_2}$. Moreover, $u_4$ is different from $u_0$, else there would be a supercycle. It is also different from $u_1$ and $u_2$, because it doesn't belong to $C_{u_2}$. Thus we have a path in $G$ with 5 vertices.\\
\indent
We can continue this process and at any step $i$, as $u_i$ is not a favourable vertex, there is a vertex $u_{i+2}$ so that $u_{i+1}u_{i+2}$ is a dominating edge and $u_{i+2}$ doesn't belong to $C_{u_i}$. We have that $u_{i+2}$ cannot be one of the $\{u_j\}_{0\leq j\leq i-2}$ because otherwise, we would have a supercycle. It is also different from $u_{i-1}$ and $u_i$, because it doesn't belong to $C_{u_{i}}$. So there is a path in $G$ with $i+2$ vertices, for any $i\in \mathbb{N}$. As the number of vertices in $G$ is finite, this process leads to a contradiction. This means that our hypothesis that there remains not covered edges after removing the last favourable vertex was wrong. Thus all the dominating edges in $G$ are covered by the cliques in $D\cup D'$.\\
\indent There are $n$ cliques in $D$. As every clique in $D'$ is
labelled by two vertices in $G$, and as no vertex in $G$ can
label two different cliques in $D'$, we deduce that $|D'|\leq
\frac{n}{2}$, proving Theorem \ref{deuxieme}. \\

\subsubsection{A bound finer than $\frac{3}{2}n$}

We remark that in the proof, we took into account the worst case when the graph induced by the edges not covered by $C=\{C^x\}_{x\in V(G)}$ had $n$ vertices. We can get a better bound if we consider the exact number of vertices of this induced graph. Let $f$ be the number of vertices in the subgraph of $G$ induced by the edges not covered by $C$. 
\begin{coroA}
If $G$ is a simple graph with $\alpha(G)=2$ containing neither supercycles nor incompatible double tops, then $ecc(G)\leq n+\frac{f}{2}$.
\end{coroA}
\begin{proof}
We have seen in the main proof that putting a clique in $D'$ implies a reduction of at least two vertices of the subgraph of $G$ induced by the edges not covered by $C$.
\end{proof}

The number $f$ can be simply obtained by substracting from $n$ the number of vertices $z$ in $G$ of which the incident edges are all in one of the following cases:
\begin{itemize}
\item it is not a dominating edge,
\item it is an edge $zx$ such that $zxy$ is a path of dominating edges and $zy$ is not an edge in $G$, $x,y\in V(G)$,
\item it is an edge $zx$ such that $xzy$ is a path of dominating edges and $xy$ is not an edge in $G$, $x,y\in V(G)$,
\item it is an edge $zx$ such that $zxy$ is a triangle of dominating edges, $x,y\in V(G)$, and there are two vertices $u,v$ in $\overline{N[y]}$ (possibly $u=v$), such that $zu$ and $xv$ are dominating edges.
\end{itemize}
This simply reflects the fact that all the edges incident to $z$ are covered by $C$.

\section{Conclusion} \label{sec:conclusion}
We proved that if a graph $G$ with $\alpha(G)=2$ has a dominating edge $uv$ and if $G\setminus\{u,v\}$ has diameter three, then the conjecture is true. We can probably go at least a step further with the same ideas and show that if $G$ has two dominating edges $uv$ and $xy$ and if $G\setminus\{u,v,x,y\}$ is of diameter three, then it also satisfies the conjecture.\\
\indent
Our second result provides a bound for graphs with no incompatible double top and no supercycles. The crucial hypothesis is that there are no incompatible double tops. It was necessary in the second step of the proof, when we put cliques in $D'$. Indeed, the fact that the not covered dominating edges that are adjacent to a given vertex form a clique supposes that every clique $C_x, x\in V(G)$ covers all the not covered edges with an end in $\overline{N[x]}$ and the other in the neighbourhood of $x$. This is only possible if the extension of $\overline{N[x]}$ to vertices $u$ and $v$ implies that $uv$ is an edge in $G$.  \\
\indent
Another, more ambitious way to find a good bound for $ecc(G)$, or even maybe the bound $n$, would be to work on maximal cliques. If the edges of a graph can be covered by $n$ cliques, then we could build on the intuition that these cliques are maximal, or can be extended to be maximal. We could try to construct step by step a set of $n$ maximal cliques that can lead to a good bound for the edge clique cover number. 

\section*{Acknowledgements}
We would like to deeply thank Ge\v{n}a Hahn and Ben Seamone for their comments that helped to improve the paper. \\
\indent
The author has been supported by the Natural Sciences and Engineering Research Council of Canada (BESC D).



\end{document}